\documentclass[10pt]{article}
\usepackage{amssymb}
\usepackage{enumerate}
\usepackage{graphicx}
\usepackage{amsthm}
\usepackage{color}
\usepackage{xcolor}

\def\R{{\rm I\! R}}

\def\C{\mbox{l\hspace{-.47em}C}}

\newtheorem{theorem}{Theorem}

\newtheorem{proposition}{Proposition}

\title{Injectivity of non-singular planar maps with one convex component}

\author{Marco Sabatini  \footnote {Dipartimento di Matematica, Univ. di Trento, I-38123 Povo (TN) - Italy; email: marco.sabatini@unitn.it. This paper has been partially supported by GNAMPA, Gruppo Nazionale per l'Analisi Matematica, la Probabilit\`a e le loro Applicazioni.}}
\date{}
\begin{document}
\maketitle
\begin{abstract}  We prove that if a non-singular planar map $\Lambda \in \C^2(\R^2,\R^2)$ has a convex component, then $\Lambda$ is injective. We do not assume strict convexity.

{\bf Keywords: } Local invertibility, global injectivity, non-strict convexity, Jacobian Conjecture. 
\end{abstract}

\section{Introduction}

Let $\Omega $ be an open connected subset of $\R^n$. We say that $\Lambda: \Omega \rightarrow \R^n$ is {\it locally injective (invertible)} at $X \in \Omega$ if there exists a neighbourhoods $U_X  \subset \Omega$ of $X$ and $V_{\Lambda(X)}$ of $\Lambda(X)$ such that the restriction $\Lambda: U_X \rightarrow V_{\Lambda(X)}$ is injective (invertible).
If $\Lambda\in C^1(\Omega,\R^n)$,  we denote by $J  (X)$ the Jacobian matrix of $\Lambda$ at $X$.
By the inverse function theorem, if  $J  (X)$ is non-singular then $\Lambda$ is locally injective at $X$. 
It is well-known that locally injective maps need not be globally injective, even if $J  (X)$ is non-singular for all $X \in \Omega$, as in the case of the exponential map $\Lambda(x,y) = (e^x \cos y, e^x \sin y)$. Injectivity (invertibility) of  locally injective (invertible) maps under suitable additional assumptions has been studied for a long time. 
 In \cite{Ke} it was conjectured that every polynomial map  $\Lambda: \C^n \rightarrow \C^n$ with  constant non-zero Jacobian determinant  be globally invertible, with polynomial inverse. Such a problem, known as {\it Jacobian Conjecture}, was widely studied and inserted in a list of relevant problems in \cite{Sm}. The Jacobian Conjecture was studied in several settings, even replacing  $\C$ with other fields, but still remains unsolved  for $n \geq 2$, \cite{BCW, BE1,dE,Ya}. In \cite{Pi} it was proved that asking for the determinant of  $J  (X)$ not to vanish is not sufficient to guarantee $\Lambda$ injectivity. 
  \\  \indent 
 Injectivity appears also in connection to a global stability problem formulated in \cite{MY}.
In such a  paper it was conjectured that if at any point $J  (X)$ has eigenvalues with negative real parts then a critical point $O$ of the differential system
\begin{equation} \label{sysn}
\dot X = \Lambda(X)
\end{equation}
 is globally asymptotically stable. Global asymptotic stability of (\ref{sysn}) implies   $\Lambda$  injectivity. In \cite{Ol} it was proved that if $n=2$, then the vice-versa is true, i. e. injectivity implies global asymptotical stability. Using such a result the conjecture was proved to be true for $n=2$ \cite{Fe,Gl,Gu}. On the other hand the conjecture does not hold in higher dimension, even for polynomial vector fields \cite{BL,CEGHM}.  
  \\  \indent 
Other additional conditions to get injectivity are growth conditions. A classical result in this field is Hadamard theorem \cite{H}, which states that if $\Lambda$ is proper, i. e. if $\Lambda^{-1}(K)$ is compact for every compact set $K \subset \R^n$, then  $\Lambda$ is a bijection. 
Properness is ensured if $\Lambda$ is norm-coercive, that is if
\begin{equation}   \label{coe}
\lim_{|X| \to +\infty} |\Lambda(X)| = +\infty.
\end{equation}
Coerciveness requires all the component of $\Lambda$ to grow enough for (\ref{coe}) to hold. On the other hand coerciveness is not necessary in order to have injectivity, as the real map $x \mapsto \arctan x$ shows. In \cite{MS1}, studying planar maps $\Lambda(z) = (P(z),Q(z))$, injectivity  was proved under a growth condition on just one component of $\Lambda$. In fact, if
\begin{equation}   \label{intgrad}
\int_0^{+\infty} \inf_{|z|=r}  |\nabla P(z)|  dr=+\infty, 
\end{equation}
then $\Lambda$ is injective. As a consequence, if there exists $k > 0$ such that $|\nabla P(z)| \geq k$, then $\Lambda$ is injective.
 \\  \indent 
Also in this paper, studying planar maps, we prove injectivity imposing a suitable condition on just one component. In fact, we prove that if one of the components $\Lambda(z) = (P(z),Q(z))$ is a non-strictly convex function, then $\Lambda(z)$ is injective. One of the steps in the proof is the same as in \cite{MS1}, since  we prove the parallelizability of the Hamiltonian system
\begin{equation}  \label{HP}
\left\{  
\matrix{\dot x = P_y \hfill  \cr  
\dot y = - P_x  }  \right.  .
\end{equation}
That is equivalent to prove the connectedness of the level sets of $P(z)$.  
 \\  \indent 
We observe that the  non-strict convexity of the function $P(z)$ implies the non-strict convexity of the orbits of (\ref{HP}), but the vice-versa is not true, as the exponential map shows. Hence injectivity cannot be proved assuming only the non-strict convexity of the orbits of 
(\ref{HP}).

\section{Maps having one convex component}  \label{convex}

In order to introduce the proof of next theorem, we recall some properties of convex functions.
\begin{proposition} \label{eleconvex} Let $f \in C^2(\R,\R)$, $H \in C^2(\R^2 , \R^2)$ be (non strictly) convex funtions. Then: 
\begin{itemize}
\item[i)] if $f$ is non-constant then it is unbounded from above;
\item[ii)] if there exist $u_1 < u_2 < u_3 \in \R$ such that $f(u_1) = f(u_2 ) = f(u_3)$, then $f$ is constant on the interval $[u_1, u_3]$;
 \item[iii)] the restriction of $H$ to every line is a convex one-variable function;
\item[iv)]  sub-level sets of $f$ and $H$ are convex;
\item[v)] every level set of $H$  at every point has a tangent line  and lies entirely on one side of such a tangent.
\item[vi)] the intersection of a level set of $H$ with any of its tangent lines is connected (a closed interval, in generalized sense).
\end{itemize}
\end{proposition}

In the proof of next theorem we consider the family of orbits of the differential system (\ref{HP}). A regular $C^1$ curve $\sigma$  is said to be a {\it section} of (\ref{HP}) if it is transversal to  (\ref{HP}) at every point of $\sigma$. 
If $\gamma$ is a non-trivial orbit, then for every $z \in \gamma$ there exists a neighbourhood $U_z$ of $z$ and  two open disjoint connected subsets $U_z^\pm \subset U_z$ lying on different sides of $\gamma$, such that $U_z = U_z^- \cup (\gamma \cap U)  \cup U_z^+$.  If $\sigma$ is a section of $\gamma$ and  $ \sigma \cap \gamma = \{ z \}$, then there exist a neighbourhood $U_z$ of $z$ and two sub-curves $\sigma^\pm$,  called {\it half-sections},  such that   $\sigma^\pm = \sigma \cap U_z^\pm$. 
 \\  \indent 
Given a planar differential system without critical points, two orbits $\gamma_1$ and $\gamma_2$ are said to be {\it inseparable} if and only if there exist two half-sections $\sigma_1$ and $\sigma_2$ such that every orbit meeting $\sigma_1$ meets also  $\sigma_2$ and vice-versa. It can be proved that if $\gamma_1$ and $\gamma_2$ are inseparable, then for every couple of points $z_1 \in \gamma_1$ and  $z_2 \in \gamma_2$ there exist half-sections such that every orbit meeting $\sigma_1$ meets also  $\sigma_2$ and vice-versa. In other words, the definition of inseparability does not depend on the choice of $z_1$ and $z_2$.
 \\  \indent 
We denote  by $\phi(t,z)$ the local flow of (\ref{HP}). Since we deal with non-singular maps, such a system has no critical points. Its orbits are positively and negatively unbounded and separate the plane into two  connected components. Every orbit is contained in a level set of $P(z)$, even if in general level sets of  $P(z)$ do not reduce to a single orbit. 
In what follows we denote by $A^o$ the interior of a set $A$ and by $\overline{A}$ its closure.

\begin{theorem} \label{teoconvex} Let $\Lambda\in C^2(\R^2,\R^2)$ be a non-singular map. If one of its components is convex, then $\Lambda$ is injective.
\end{theorem}
{\it Proof.} 
Possibly exchanging the components, we may assume $P(z) $ to be convex. By lemma 2.2 and theorem 2.1 in  \cite{MS1}, it is sufficient to prove that the level sets of $P(z)$ are connected.
By absurd, let us assume that a level set of $P(z) = h$ is disconnected. As a consequence by lemma 2.2 in \cite{MS1} the system (\ref{HP}) has a couple  $\gamma_1 \neq \gamma_2$ of inseparable orbits. By continuity, $P(z)$ assumes the same value on $\gamma_1 $ and 
$\gamma_2$, say  $P(\gamma_1) = P(\gamma_2) = k$. 
 \\  \indent 
Let us consider two cases.
 \\  \indent 
1) One among $\gamma_1 $ and $ \gamma_2$ is not a line. Assume $\gamma_1 $ is not a line.  Let $\Gamma_1$ be the closed convex set having $\gamma_1$ as boundary.
 \\  \indent 
1.1)  If $\gamma_2 \subset \Gamma_1$, then it is not a line, otherwise it would meet $\gamma_1$, contradicting uniqueness of solutions. Let $z_1$ be an arbitrary point of  $\gamma_1$ and $\tau_{12}$ be the line passing through $z_1$ and tangent to $\gamma_2$, existing by the convexity of $\Gamma_2$. Since $\gamma_2$ is not a line one can rotate  $\tau_{12}$ around $z_1$ until it meets $\gamma_2$ at two points $z_2^1 \neq z_2^2$. Let us call $\tau^*$ such a line. Then $\tau^*$ meets the level set  $P(z) = k$ at three distinct points, $z_1, z_2^1, z_2^2$.  By proposition \ref{eleconvex}, $ii)$, $P(z)$ is constant on the smallest segment $\Sigma$ containing $z_1, z_2^1, z_2^2$. The set $\gamma_1 \cup \Sigma \cup \gamma_2$ is connected and contained in $P(z) = k$, contradicting the fact that $\gamma_1 $ and $ \gamma_2$ are distinct connected components of $P(z) = k$. 
 \\  \indent 
1.2) Let $\gamma_2 \subset \Gamma_1^c$. If $\gamma_1 \subset \Gamma_2$, then one can reply the argument of point 1.2), exchanging the role of $\gamma_1 $ and $\gamma_2$.
 \\  \indent 
 1.3) Assume $\gamma_1 \not\subset \Gamma_2$ and $\gamma_2 \not\subset \Gamma_1$.  
Let $D_1$ be the subset of $\gamma_1$ consisting of its linear parts, i.e. half-lines and line segments. Since $\gamma_1 $ is not a line, one has $D_1 \neq \gamma_1$. Let us choose arbitrarily $z_1 \in \gamma_1 \setminus D_1$ and let $\tau_1$ be the  tangent line of $\gamma_1$ at $z_1$. By  point $v)$ of Proposition \ref{eleconvex} $\gamma_1$ lies on one side of $\tau_1$. One has $\gamma_1 \cap \tau_1 = \{ z_1 \}$.  Let $\tau_1^\pm$ be the half-lines contained in $\tau_1$ having $z_1$ as extreme point, $\tau_1^+$ tangent to the positive semi-orbit of $z_1$,  $\tau_1^-$ tangent to the negative semi-orbit of $z_1$. 
Let $\Pi_1$ the closed half-plane having $\tau_1$ as boundary and containing $\gamma_1$. For all $\epsilon > 0$ one has $\phi(\pm \epsilon, z_1) \in \Pi_1^o$. Every such orbit meets $\tau_1$ at least at two points lying on distinct half-lines. As a consequence, $z_1$ is an isolated point of minimum of the restriction of $P(z)$ to the line $\tau_1$. Hence $\gamma_2$ does not meet $\tau_1$.  \\
By the inseparability of $\gamma_1 $ and $ \gamma_2$ there are half-sections $\sigma_1$ of $\gamma_1 $ at $z_1$ and $\sigma_2$ of $ \gamma_2$ at $z_2$ such that every orbit meeting $\sigma_1$ meets also $\sigma_2$ and vice-versa. One can take $\sigma_1$ and $\sigma_2$ small enough to have $\overline{\sigma_1}$ and  $\overline{\sigma_2}$ compact, disjoint and such that   $\sigma_2 \cap \Pi_1 = \emptyset$.
 \\  \indent 
There exist neighbourhoods $U_\epsilon^\pm$ of $\gamma_1(\pm \epsilon) $ such that $U_\epsilon^\pm \subset\Pi^o$. 
By the continuous dependance on initial data there exists a neighbourhood $U_1$ of $z_1$ such that $\phi( \pm \epsilon ,U_1) \subset U_\epsilon^\pm$. This holds in particular for the points of $\delta_1 = \sigma_1 \cap U_1$, so that  $\phi( \pm \epsilon ,\delta_1) \subset U_\epsilon^\pm \subset \Pi_1^o$. $\delta_1$ is itself a half-section at $z_1$.
For all $z \in \delta_1$ the orbit $\phi(t,z)$ meets both $\tau_1^-$ and $\tau_1^+$, hence both half-lines contain points $z^\pm$ such that  $P(z^-) = P(z^+) > P(z_1)$. Moreover, $\phi(t,z)$ does not meet $\tau_1$ at a third point, since in that case, by point $ii)$ of Proposition \ref{eleconvex}, $P(z)$ would be constant on a segment of $\tau_1$ containing $z_1$, contradiction. Hence, for all $z \in \delta_1$, both semi-orbits starting at $z$ are definitively (resp. for $t \to \pm \infty$) contained in $\Pi_1^o$. \\ \indent
The set $W = \phi( [ - \epsilon,  \epsilon], \overline{\delta_1}) $ is compact. It is possible to take $\overline{\delta_1}$ small enough in order to have $z_2 \not \in W \cup \Pi_1$ (otherwise $z_2 = z_1$).
 By construction, every orbit starting at a point of $\overline{\delta_1}$ is contained in the closed set $  W \cup \Pi_1$. Let us denote by $\delta_2$ the part of $\sigma_2$ met by orbits starting at points of $\delta_1$. Since every point of $\delta_2$ lies on an orbit starting at $\delta_1$, the half-section $\delta_2$  is contained in  $W \cup \Pi_1$. As a consequence, one has
$$
z_2 \in \overline{\delta_2} \subset W \cup \Pi,
$$
contradiction.
  \\  \indent 
 2) Assume both $\gamma_1$ and $\gamma_2$ to be  lines. They are parallel, since otherwise they should meet at a point $z_0$ which should be a fixed point of (\ref{HP}), contradicting the nonsingularity of $\Lambda$. Let $\Sigma_{12}$ be the closed strip having boundary $\gamma_1 \cup \gamma_2$. 
Let $\sigma $ be a line orthogonal to $\gamma_1$ and $\gamma_2$, and let us set $z_1 = \gamma_1 \cap \sigma$, $z_2 = \gamma_2 \cap \sigma$,  $\sigma_{12} =  \Sigma_{12} \cap \sigma$. The orbits   $\gamma_1$ and $\gamma_2$ are inseparable, hence there exist open sub-segments $\sigma_1$ and $\sigma_2$ of  $\sigma_{12}$ such that 
$z_1 \in \overline{\sigma_1}$, $z_2 \in \overline{\sigma_2}$,  $\overline{\sigma_1} \cap \overline{\sigma_2} = \emptyset$ and every orbit meeting $\sigma_1$ meets $\sigma_2$, and vice-versa. 
Let $\Phi_{12}$ be the union of the orbits meeting $\sigma_1$ and $\sigma_2$. Both $\gamma_1$ and $\gamma_2$ are contained in $\partial\, \Phi_{12}$. The restriction of $P(z)$ to the compact set $\sigma_{12}$ is convex and non constant (because if it was constant $\gamma_1$, $\gamma_2$ and $\sigma_{12}$ would be in $P(z) = k$, contradiction). One has 
$$
 \max \{ P(z) : z \in \sigma_{12} \} =P(z_1) = P(z_2) = k .
$$
Let $z_m$ a point of $\sigma_{12}$ such that 
$$
P(z_m) = \min \{ P(z) : z \in \sigma_{12} \} < P(z_1) = P(z_2) = k.
$$
The orbit starting at $z_m$ is tangent to $\sigma_{12}$ and lies entirely on one side of $\sigma_{12}$. One has $\nabla P(z_m)\, \bot\, \sigma_{12}$, with the vector $\nabla P(z_m)$ pointing towards the half-strip $\Sigma_{12}^+$ not containing $\phi(t,z_m)$. 
Let $\eta$ be the line  parallel to $\gamma_1$ and $\gamma_2$ passing through $z_m$. The line $\eta$ meets all the orbits passing through $ {\sigma_1}$ and ${\sigma_2}$, hence the restriction of $P(z)$ to $\eta$ assumes every value belonging to $[P(z_m),k)$. 
On the other hand, by proposition \ref{eleconvex}, $i)$,  $P(z)$ is unbounded from above on $\eta$, hence there exists a point in $z\in\eta$ such that $P(z)= k$. Let $z_{12}$ the point such that $P(z_{12})= k$, closest to $z_m$. Then the orbit $\phi(t,z_{12})$ is inseparable from 
$\gamma_1$ and $\gamma_2$, since every orbit meeting $ {\sigma_1}$ and ${\sigma_2}$ also meets $\eta$ in a neighbourhood of $z_{12}$. In other words, a suitable sub-segment $\eta_{12}$ of $\eta$ is a half-section of $\phi(t,z_{12})$ such that every orbit meeting $\sigma_1$ and $\sigma_2$ meets also $\eta_{12}$, and vice-versa. 
 \\  \indent 
The orbit $\phi(t,z_{12})$  cannot be a line because in such a case either it would be parallel to $\gamma_1 $ and $ \gamma_2$, contradicting their inseparability, or transversal to them, implying the existence of two critical points, $\gamma_1 \cap \gamma_{12}$ and $\gamma_2 \cap \gamma_{12}$. Since $\gamma_{12}$ is not a line point 1) applies.
  
 \hfill  $\clubsuit$

A simple example of non-linear non-singular map with both non-strictly convex components is
$$
\Lambda(x,) = (x+y+e^x,x+y+e^y).
$$
The Hamiltonian system of a non-strictly convex two-variables function has non-strictly convex orbits. The vice-versa is not true, as the function $e^x \cos y$ shows. Infact, the connected components of  $e^x \cos y = 0$ are lines, and the connected components of  $e^x \cos y = k \neq 0$ are strictly convex, since they are graphs of the one-variable functions
$$
x = \ln \left(   \frac k{\cos y}   \right),
$$
whose second derivative does not vanish. On the other hand the hessian matrix of $e^x \cos y$ is:
$$
\left ( 
\matrix{ e^x \cos y & - e^x \sin y \hfill  \cr  
 - e^x \sin y & - e^x \cos y  }  \right )  ,
$$
whose Jacobian determinant is $- e^{2x} < 0$. In fact, the map $\Lambda(x,y) = (e^x \cos y,$ $ e^x \sin y)$ is not injective, even if both Hamiltonian systems of its components have non-strictly convex orbits.

\end{document}